\documentclass[11pt]{amsart}

\usepackage{amssymb}
\usepackage{graphics}
\usepackage{graphicx}
\usepackage{amscd}

\newtheorem{theorem}{Theorem}
\newtheorem{lemma}{Lemma}
\newtheorem{prop}{Proposition}

\newtheorem{cor}{Corollary}
\newtheorem{conj}{Conjecture}

\theoremstyle{definition}
\newtheorem{remark}{Remark}
\newtheorem{example}{Example}

\newcommand{\on}{\operatorname}

\newcommand{\Spinc}{\on{Spin}^c}
\newcommand{\gr}{\on{gr}}

\newcommand{\h}{\widehat}
\newcommand{\CKh}{\widetilde{CKh}}
\newcommand{\Kh}{\widetilde{Kh}}
\newcommand{\uu}{{\bf u}}
\newcommand{\um}{{\bf u}_-}
\newcommand{\up}{{\bf u}_+}

\newcommand{\OO}{\mathcal{O}}
\newcommand{\f}{\frac}
\newcommand{\Z}{\mathbb{Z}}
\newcommand{\R}{\mathbb{R}}
\newcommand{\Q}{\mathbb{Q}}
\newcommand\goth[1]{\mathfrak{#1}}
\newcommand{\s}{\goth{s}}

\newcommand{\sign}{\on{sign}}

\newcommand{\unknot}{\mbox{unknot}}
\newcommand{\ps}{\ti{\psi}}
\newcommand{\ti}{\tilde}

\begin{document}

\author{Olga Plamenevskaya}
\address{Department of Mathematics, Stony Brook University, Stony Brook, NY 11794}
\email{olga@math.sunysb.edu}
\title[Transverse knots and Heegaard Floer contact invariants]{Transverse knots, branched double covers and Heegaard Floer contact invariants}

\begin{abstract}
Given a transverse link in $(S^3,\xi_{std})$, we study the contact
manifold that arises as a  branched double cover of the sphere. We
give a contact surgery description of such manifolds, which allows
to determine the Heegaard Floer contact invariants for some of
them. By example of the knots of Birman--Menasco, we show that
these contact manifolds may fail to distinguish between non-isotopic
transverse knots. We also investigate the relation between the
Heegaard Floer contact invariants of the branched double covers
and the Khovanov homology, in particular, the transverse link
invariant we introduce in a related paper.

\end{abstract}

\maketitle

\section{Introduction}
Let $(S^3, \xi_{std})$ be the 3-sphere equipped with its
standard contact structure $\xi_{std} =\ker(dz-xdy)$. A link
$L\subset S^3$ is called {\em transverse} if it is everywhere
transverse to the contact planes. Let $\Sigma(L)$ be the double
cover of $S^3$ branched over $L$. Then, $\Sigma(L)$ carries a
natural contact structure $\xi_L$ lifted from $(S^3,\xi_{std})$.
The goal of this paper is to study the contact manifold
$(\Sigma(L), \xi_L)$.

Our motivation is two-fold. First, one might wonder whether the
double branched covers can help us understand transverse knots and
links. The classification of transverse knots is a very difficult
task. Indeed, while a few simplest knots, such as the unknot, the
figure eight knot, and torus knots are completely classified
by their topological knot type and the self-linking number \cite{EH1}, in
general this is not true. The first examples of smoothly isotopic,
but not transversely isotopic transverse knots $K_1$, $K_2$ with
the same self-linking number were given by Birman and Menasco
\cite{BM}. The existence of such pairs was also demonstrated  by
Etnyre and Honda \cite{EH3} via some ``non-explicit" examples.
Unfortunately, it seems that the branched double covers do not
capture the subtle difference between such knots. Indeed, we prove

\begin{theorem} \label{bm1} Let $K_1$, $K_2$ be transversely non-isotopic
knots of \cite[Theorem 3]{BM}. The branched double covers
$(\Sigma(K_1), \xi_{K_1})$ and $(\Sigma(K_2), \xi_{K_2})$ are
contactomorphic.
\end{theorem}
(Theorem 3 in \cite{BM} provides a family of pairs of transversely
non-isotopic knots, not just a single pair. Our result is true for
all pairs of  \cite{BM}).

From another viewpoint, double covers of $(S^3, \xi_{std})$ branched over transverse
links give an interesting
special case of contact 3-manifolds.
Indeed, the works of Giroux  \cite{Gi} and others imply that every 3-dimensional
contact manifold can be represented as a triple branched cover of
$(S^3, \xi_{std})$. It turns out that  branched double covers are quite simple.
We give an algorithm for finding a contact surgery diagram
\cite{DG,DGS} for $(\Sigma(L), \xi_L)$; the diagrams we get only
involve surgeries on Legendrian unknots. We are also able to find
the homotopy invariants of the contact structure $\xi_L$ (that is,
the induced $\Spinc$ structure and the three-dimensional invariant
$d_3(\xi)$ of \cite{Go}).
\begin{theorem} Let $\s_L$ be the $\Spinc$ structure induced by
$\xi_L$. Then $c_1(\s_L)=0$. The  invariant $d_3(\xi)$ is
completely determined by the topological link type of $L$ and its
self-linking number $sl(L)$.
\end{theorem}
In certain cases, it is easy to tell whether the contact structure $\xi_L$ is tight or
overtwisted.
\begin{prop} \label{qua-pos}
 $(\Sigma(L), \xi_{L})$ is overtwisted if $L$ is obtained as a
transverse stabilization of another transverse link.
\end{prop}

\begin{prop} \label{stab} $(\Sigma(L),\xi_L)$ is Stein fillable if the transverse link
$L$ is represented by a quasipositive braid.
\end{prop}
(Here and later on, it will be convenient to represent transverse
links as closed transverse braids in $(S^3, \xi_{std})$. We give
more details about this representation
 in the next section.)

We now turn attention to the Heegaard Floer contact invariants.
The Heegaard Floer theory of Ozsv\'ath and Szab\'o (\cite{OS1} and sequels) associates a
homology group $\widehat{HF}(Y)$ to a closed oriented 3-manifold
$Y$ and yields invariants for many low-dimensional objects. In
particular, given a contact structure $\xi$ on $Y$, the contact
invariant $c(\xi)$ is a distinguished element of
$\widehat{HF}(-Y)$, defined up to sign \cite{ContOS}. (We assume that the
coefficients are taken in  $\Z$.)

Propositions \ref{qua-pos} and \ref{stab} along with the surgery
diagrams  and properties of  $c(\xi)$ enable us to determine the
Ozsv\'ath-Szab\'o contact invariant $c(\xi_L)$ for
$(\Sigma(L),\xi_L)$ in many cases. We observe that the contact
invariant $c(\xi_L)$ behaves very similarly to the
Khovanov-homological invariant of transverse links that we
introduce in \cite{Pla}. This is not a mere coincidence. Indeed,
as proved by Ozsv\'ath and Szab\'o \cite{OSKh}, for a smooth link
$L\subset S^3$ there is a spectral sequence converging to
$\widehat{HF}(-\Sigma(L))$ whose $E^2$ term is given by the reduced
Khovanov homology of $L$ (both theories are to be taken with
$\Z/2\Z$ coefficients). When the link $L$ is alternating, the
spectral sequence of \cite{OSKh} collapses at the $E^2$ stage,
providing an isomorphism 
\begin{equation} \label{isom}
\widetilde{Kh}(L) \cong \widehat{HF}(-\Sigma(L)).
\end{equation}
Now, suppose the link $L$ is transverse, and $\xi_L$ is the induced 
contact structure on the branched double cover. Let  $\psi(L)\in
\widetilde{Kh}(L)$ be the invariant of \cite{Pla}.
We would like to suggest that the elements $\psi(L)$ and   $c(\xi_L)$
correspond to one another under the isomorphism (\ref{isom}). However, we must be careful, 
because this isomorphism is not canonical:   
while the spectral sequence of \cite{OSKh} is believed to be an invariant of the link,
this invariance has not been proved.  To deal with this issue, 
we fix a link diagram before studying such an isomorphism. 
(The choice of the diagram  will be clear from the context; besides, we prove that 
both $\psi(L)$ and  $c(\xi_L)$ are independent of the diagram.) 

We also need to be more precise about the spectral sequence and  isomorphism~(\ref{isom}).
For a fixed link diagram, the construction of \cite{OSKh} gives a filtered chain complex $C(L)$, whose homology 
is $\h{HF}(-\Sigma(L))$, and the associated graded complex is the chain complex $\widetilde{CKh}(L)$ for reduced 
Khovanov homology (with its homological grading). When $L$ is alternating, the  spectral sequence collapses, 
yielding  a canonical isomorphism between $\widetilde{Kh}(L)$ and the associated graded
group of $\h{HF}(-\Sigma(L))$. For coefficients in a field, the associated graded
group of $\h{HF}(-\Sigma(L))$ is of course isomorphic to $\h{HF}(-\Sigma(L))$, but the latter isomorphism is 
not canonical. Therefore, we will always think of~(\ref{isom}) as the isomorphism between $\widetilde{Kh}(L)$ and the associated graded group of $\h{HF}(-\Sigma(L))$. The relation between $\psi(L)$ and  $c(\xi_L)$ must 
involve gradings, as follows. Recall from \cite{Pla} that $\psi(L)$ is a homogeneous element of $\widetilde{Kh}(L)$
of homological degree $0$. Let $c_0 (\xi_L)$ be the image of $c(\xi_L)$ in the corresponding 
subquotient of the associated graded group of  $\h{HF}(-\Sigma(L))$, that is, 
$c_0 (\xi_L)\in \h{HF}_0(-\Sigma(L))/ \h{HF}_1(-\Sigma(L))$, where 
the subscripts on $\h{HF}$ indicate the filtration level.

We suggest
\begin{conj} \label{conjec} If $L$ is a transverse representative of an
alternating smooth link, then the homological grading of 
$\psi(L) \in \widetilde{Kh}(L)$  is the same as the filtration level 
of $c(\xi_L) \in  \h{HF}(-\Sigma(L))$ (and this filtration level is $0$).
Moreover, $\psi(L)=  c_0 (\xi_L)$ under the isomorphism (\ref{isom})
between  $\widetilde{Kh}(L)$ and the associated graded group of 
$\h{HF}(-\Sigma(L))$.
\end{conj}

In the general case, it is plausible that $c(\xi_L)$ somehow
``corresponds" to $\psi(L)$ under the spectral sequence.

In the special case when the transverse link $L$ is
represented by a transverse closed braid whose braid diagram is
alternating,  Conjecture \ref{conjec} is not hard to prove.
\begin{theorem} \label{psi=c} Let $L$ be a transverse link represented by
a closed braid with  an alternating braid diagram. Then, the filtration level of $c(\xi)$
is as stated, and $\psi(L)=  c_0 (\xi_L)$. 
\end{theorem}

It should be noted that alternating braids represent a very narrow
class of links. We show that $\psi(L)=c(\xi_L)=0$ for all such
links except the $(2, n)$-torus link.

Conjecture \ref{conjec} implies that $c(\xi_L)\neq 0$ whenever $\psi(L) \neq 0$.
We are able to check this directly for many cases not covered by Theorem \ref{psi=c}
(see section \ref{OSi}). If proved in general,  Conjecture \ref{conjec}, together 
with the fact that $c(\xi)$ vanishes when $\xi$ is overtwisted \cite{ContOS}, would give 
a powerful sufficient condition for the contact structure $\xi_L$ to be tight.
Indeed, for an arbitrary transverse link $L$ we can often show that $\psi(L) \neq 0$
by using arguments from \cite{Pla} or software \cite{BN}, \cite{Sh}; when  $K$ is an 
alternating knot, it can be shown that $\psi(K)\neq 0$ if and 
only if $sl(K) =-\sigma(K) -1$, where $\sigma(K)$ is the signature of the knot 
(with the sign convention such that the right-handed trefoil has signature $-2$.)
Note that the contact structures that we would thus obtain from  Theorem \ref{psi=c} 
are all (trivially) Stein fillable.

\noindent {\bf Acknowledgements.} I would like to thank Peter
Ozsv\'ath 
and Andras Stipsicz  for very helpful discussions. I am also very grateful for the referee's 
remarks and suggestions.

\section{Transverse links and braids}

In what follows, we will be working with 
the induced contact structure $\xi_L$ on the branched double cover $\Sigma(L)$
for a transverse link $L$ in the standard contact 3-sphere.
We now  describe in some detail how $\xi_L$ is constructed.

Let $L$ be a transverse knot in $S^3$ (if $L$ is a link, we can deal with every component
separately). Then some neighborhood of $L$ embeds into $\R^2\times S^1$ 
via coordinates $(r, \theta, z)$
(with $(r, \theta)$ the polar coordinates on $\R^2$, $z\in S^1$, and $L= \{r=0\}$), and 
 the contact structure $\xi$ in this neighborhood can be given as the kernel of 
the 1-form $dz+r^2d\theta$. (This is the Darboux theorem for contact structures.) In this neighborhood, 
the standard local model for the projection $p:\Sigma(L) \to S^3$ is given by the map $(w, z) \to (w^2, z)$, 
where $w=x+iy$. 

Fix  $n$ large enough, so that the set  $\{ r = r_0\}$ with $r_0^4=\frac1{4 \pi n}$  is contained in the chosen neighborhood. 
Choose $\eta, \epsilon >0$ such that $\eta< r_0$, and $\epsilon^2< 2 \eta^4$. Define $\xi_L$ on $\Sigma(L)-\{r<\eta\}$ 
to be the kernel of the pull-back contact form (i.e. $\xi_L= \ker(dz+2r^4 d\theta)$ where our coordinates are defined, 
and $r\geq \eta$). For $r<\epsilon$, let $\xi_L$ be the kernel of the contact form $dz+r^2 d\theta$, and interpolate 
between the two pieces by setting $\xi_L=\ker (dz+ f(r)d\theta)$, where the smooth function $f$ is chosen so that 
$f(r)=r^2$ for $r<\epsilon$, $f(r)=2r^4$ for $r>\eta$, and $f'(r)>0$ for $r>0$. It is clear that $dz+ f(r)d\theta$
is a contact form; moreover, the contact structure it defines inside the coordinate neighborhood of $L$ in $\Sigma(L)$
is isotopic to $\ker(dz+r^2 d\theta)$, and therefore tight.  (Note that the pull-back form $dz+2r^4 d\theta$ would not 
work for the entire $\Sigma(L)$: this form is not contact along the $z$-axis!)

We have to check that the contact structure $\xi_L$ we obtain is independent of choices.
To this end, observe that the characteristic foliation that $\xi_L$ induces on  the torus $\{ r = r_0\}$ 
is given by parallel longitudes of framing $-n$ (calculated with respect to the framing defined by $\theta$). 
Fix two of these parallel longitudes. They divide the torus $\{ r = r_0\}$ into two annular regions; by pushing one of these 
regions in and another one out, we can perturb this torus by an isotopy into a convex surface whose 
characteristic foliation is Morse--Smale, and the dividing set is given by the two parallel 
longitudes of framing $-n$. In addition, we can assume that the support of this  isotopy lies outside the set $\{r\leq \eta\}$.   
 We denote the new convex torus by $T_n$, and the tubular 
neighborhood of $L$ in $\Sigma(L)$ that $T_n$ bounds by $V_n$.  

With given boundary conditions, the tight contact structure in the solid torus bounded by $T_n$
is unique up to an isotopy (see e.g. \cite{Ho}). This means  that the  tight contact 
structure on $V_n$  is uniquely determined by the boundary conditions,
and, since our construction of $\xi_L$ outside of $T_n$ is canonical, the contact structure $\xi_L$ on $\Sigma(L)$
is independent of all choices (including the choice of $n$).

It will be helpful to think about transverse links by representing them 
by closed braids. For this, consider the symmetric version of
$(S^3,\xi_{sym})$ with $\xi_{sym}=\ker(dz+xdy-ydx)$. Then, any
closed braid around $z$-axis can be made transverse to the contact
planes; moreover, any transverse link in $(S^3, \xi_{std})$ is
transversely isotopic to a closed braid \cite{Be}.

To define the self-linking number $sl(L)$, trivialize the plane field $\xi$,
and let the link $L'$ be the push-off of $L$ in the direction of the first coordinate
vector for $\xi$. Then, $sl(L)$ is the linking number between $L$ and $L'$.
Given a closed braid representation of $L$, we have
\begin{equation}\label{sl-braid}
sl(L)=n_+ -n_- - b,
\end{equation}
where $n_+$ ($n_-$) is the number of positive (negative) crossings,
and $b$ is the number of strings in the braid.

The {\em stabilization} of a transverse link represented as a braid is
 equivalent to the negative braid stabilization, i.e. adding an extra string
and a negative kink to the braid.
If $L_{stab}$ be the result of
stabilization of $L$, then
\begin{equation}
sl(L_{stab})=sl(L)-2.
\end{equation}
Note that the positive braid stabilization does not change the transverse type of the link.

We will often describe a braid by its braid word on the standard
generators $\sigma_1, \sigma_2,\dots$ and their inverses, and draw associated braid
diagrams.






\section{Crossing Resolution and Contact Surgery} \label{surgery}

In this section, we establish the correspondence between the
crossings of the braid diagram representing the transverse link
$L$ and the surgeries required to obtain the contact manifold
$(\Sigma(L), \xi_{L})$ from $(S^3, \xi_{std})$. This
correspondence is our main tool: it allows to build
contact surgery diagrams and plays the key role in other results
of this paper.

\subsection{Contact surgery}
Let $K$ be a null-homologous Legendrian knot in a contact manifold
$(Y, \xi)$. Legendrian surgery on $K$ is the surgery with
coefficient $tb(K)-1$; it is well-known that a Legendrian surgery
on  $(Y, \xi)$ produces a new contact manifold $(Y', \xi')$, which
is Stein  fillable if $(Y, \xi)$ is. We often refer to Legendrian
surgery as $(-1)$ contact surgery (comparing the surgery framing
to the framing given by contact planes). It is also possible to
make sense of $(+1)$ contact surgery (in fact, any rational $p/q$
surgery) \cite{DG}. The new contact manifold is obtained by
cutting out a tubular neighborhood of the knot $K$ (i.~e. a solid
torus) and gluing it back in so that the contact structure on the
solid torus matches the contact structure on its complement; when
the surgery coefficient is $1/q$ with $q\in \Z$, the result of
this procedure is independent of choices. We refer the reader to
\cite{DG,DGS} for the details of this construction. We recall that
$(+1)$ contact surgery is the operation inverse to the $(-1)$
contact (i.e. Legendrian) surgery, and note that the $(+1)$
surgery does not preserve Stein fillability or other similar
properties of contact structures.

\subsection{Surgery diagrams for double covers}
We are now ready to relate the crossings of the braid diagram
 to contact surgeries. Roughly, a positive crossing gives a $(-1)$-surgery,
 and a negative crossing a $(+1)$-surgery.

\begin{theorem} \label{resolve}
{\em (1)}  If the transverse braid
$L$ is obtained  from the transverse braid $L_+$  by resolving a positive
 crossing, then $(\Sigma(L_+),\xi_{L_+})$ is obtained from $(\Sigma(L),\xi_L)$
by Legendrian surgery.

{\em (2)} If $L$ is obtained from $L_-$ by resolving a negative crossing, then
the contact manifold  $(\Sigma(L_-),\xi_{L_-})$ is obtained from $(\Sigma(L),\xi_L)$
by $(+1)$ contact surgery.
\end{theorem}

\begin{proof}

We first 
 consider a  model example of two simple braids in $(S^3, \xi_{std})$.
Let $K\subset S^3$ be the transverse unknot given by the braid $\sigma_1$, and $K_+$ be the transverse
Hopf link given by $\sigma_1^2$. 
We claim that  $\Sigma(K)=S^3$ and $\Sigma(K_+)=L(2,1)=\R P^2$ with their (unique)
tight contact structures. Indeed, since $K$ is the transverse unknot with $sl(K)=-1$, it can be thought of as the binding of an open book
decomposition of $S^3$ whose page is a disk. The branched double cover, then, is the same open book,
giving the standard contact structure on $S^3$. The positive Hopf link $K_+$ with $sl(K_+)=0$  is
the binding of an open book whose page is an annulus, and the monodromy is the positive Dehn twist.
Doubling the monodromy, we see that $\Sigma(K_+)=L(2,1)$, and that the contact structure $\xi_+$
is Stein fillable. (See \cite{Gi} for the relation between open book decompositions and
fillability.)

Now,  $K$ and $K_+$ agree in the complement of a  ball $B$ containing one of the crossings
of $K_+$ (the boundary of this ball is shown in Fig. \ref{two-crossings} as a dotted circle).
The double cover of this complement $S^3\setminus B$ branched along the two arcs of $K$
is a solid torus $M$.
\begin{figure}[ht]
\includegraphics[scale=0.7]{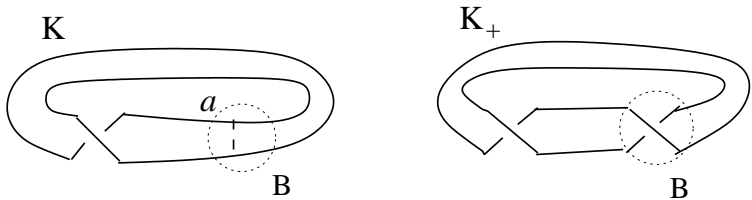}
\caption{The branched double covers for the unknot and the Hopf link differ by regluing a solid torus which is a double cover of the ball $B$. This amounts
to a surgery on the knot obtained as the branched double cover of the arc $a$.}
\label{two-crossings}
\end{figure}
The ball $B$ contains two arcs of $K$ and two arcs of $K'$, and the branched double
covers $\Sigma(K)$ and $\Sigma(K_+)$ are both obtained from $M$ by attaching another solid torus $N$,
a double cover of $B$ branched along two arcs. This means that $\Sigma(K_+)$ is obtained
as a surgery on a knot in $\Sigma(K)$. This knot is the branched double cover $\ti a$ of the arc $a$
connecting the two strands of $K$ inside $B$; it represents the longtitude of the solid torus $N$.
Therefore, we are doing surgery on an unknot in $S^3$; since the result of this surgery is $L(2, 1)$,
the surgery coefficient is $-2$. (For surgeries relating branched double covers, the surgery coefficient is always 
integral.)

Our goal is to put this surgery into a contact context. First, we can  assume that the arc $a$
is Legendrian (and so is $\ti a$). The surgery on the contact $S^3$ then becomes contact surgery on 
a Legendrian unknot; if $tb(\ti a)= -n$ for some $n\geq 1$, then we must be doing $(n-2)$-contact surgery 
(the result of this surgery is not unique unless $n-2 =\pm 1$). Our goal is to show that $n=1$, 
so that we have the Legendrian surgery on the standard Legendrian unknot in $S^3$.
To this end, we will show that the resulting contact structure is overtwisted whenever $n>1$.
Indeed, consider  $(+1)$-contact surgery on an unknot $\ti a$ with $tb(\ti a)=-3$ (Fig. \ref{tb3}). 
\begin{figure}[ht]
\includegraphics[scale=0.65]{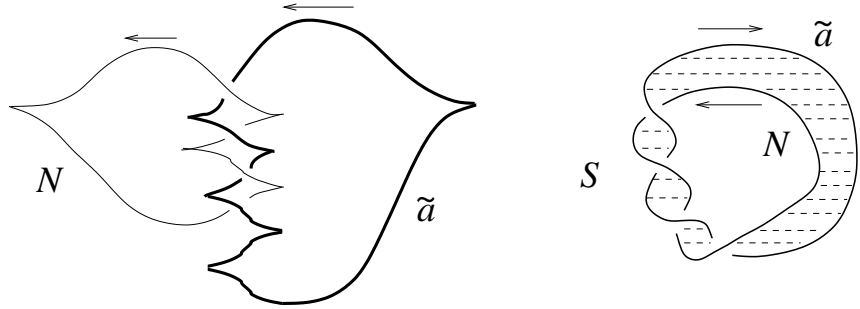}
\caption{The unknot $N$ bounds an overtwisted disk in the surgered manifold; part of the disk is 
formed by the surface $S$.}
\label{tb3}
\end{figure}
The surface framing of  $\ti a$ given by the Seifert surface $S$ of the link $N\cup \ti a$ 
in this figure equals  $-2$, which is the Dehn surgery coefficient. Then the new meridional disk 
 glued together with $S$ gives an embedded disk $D_0$ bounded by $N$ in the surgered manifold.
The  Thurston--Bennequin number of $N$ is $-2$, which is the same 
as  the surface framing of $N$ determined by the disk $D_0$.  It follows that $D_0$ is an overtwisted disk.
(A similar argument for an overtwisted sphere is given in \cite{DGS}).  
There are two more unknots with $tb=-3$, one with three kinks on the right and one on the left, and 
the other with two kinks on each side; they can be treated in the same way (putting $N$ on the left 
or on the right, and orienting the link as needed).
\begin{figure}[ht]
\includegraphics[scale=0.65]{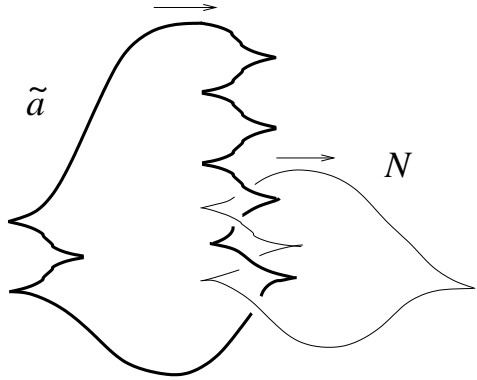}
\caption{As in the previous picture, the unknot $N$ bounds an overtwisted disk in the surgered manifold.}
\label{tbn}
\end{figure}

If  $tb(\ti a)=-n<-3$, the $(n-2)$-contact surgery is no longer uniquely defined, but the same   
argument applies, and we can aways find an unknot $N$ that bounds an overtwisted disk in the surgered manifold.
One possible picture for $\ti a$ and $N$ is shown of Fig. \ref{tbn} (note that there a few choices 
for the contact surgery even when $\ti a$ is fixed, see \cite{DGS}).    
The remaining case is $tb(\ti a)=-2$. Then, the surgery coefficient matches the Thurston--Bennequin number
(i.e. we are attempting to do $0$-contact surgery); in this case, our Legendrian unknot bounds an overtwisted 
disk in the surgered manifold.
    
Therefore, we see that the contact manifold $(\Sigma(K_+),\xi_+)$ is the result of the Legendrian
surgery on an unknot  in $(\Sigma(K),\xi)$ with $tb=-1$. This unknot  is a Legendrian
representative of the lifting of a chord $a$ of $K$, shown in  Fig. \ref{two-crossings}.
 Because (+1)-contact surgery and Legendrian surgery are inverse to one another, 
we can also say that $(\Sigma(K),\xi)$ is obtained from $(\Sigma(K_+),\xi_+)$ by the contact (+1)-surgery on $\ti a$.

We are now ready to prove part (1) of the theorem. We know that the contact manifold 
$(\Sigma(L_+),\xi_{L_+})$ is obtained by some contact surgery on the Legendrian knot $\ti a$
 obtained as the branched double cover of an  arc connecting two strings of the braid $L$.
The framing of this surgery is a purely local question, so from our model example, we see
that it must be the Legendrian framing.

Part (2) also follows from the local model: 
indeed, the link $K_+$ can be obtained from $K$
by resolving a negative crossing (we first change the diagram by a Reidemeister move
to introduce two new crossings, a positive and a negative one, that cancel each other).
\end{proof}



We can now give a contact surgery diagram for a double cover branched over an arbitrary
 transverse braid  $L$  on $n$ strings.
Inserting the factors $\sigma_i\sigma_i^{-1}$ if necessary, we may assume that the braid word
contains each $\sigma_i$ with $i=1,\dots,n$. Therefore, the braid $L$ can be obtained from
the braid $U=\sigma_1\sigma_2\dots\sigma_n$ by introducing some extra (positive and negative) crossings.
Obviously, the transverse link $U$ is the unknot with $sl(U)=-1$,
and so the double cover of $S^3$ branched over $U$ is $(S^3, \xi_{std})$.
Now, $(\Sigma(L),\xi_L)$ is the result of a
contact surgery on a link in the standard sphere:
each ``extra" crossing of $L$ gives a $(-1)$ or a $(+1)$ contact surgery on a Legendrian unknot
(depending on the sign of the crossing). These Legendrian unknots are the components for the surgery link
for $(\Sigma(L),\xi_L)$; it remains to understand how they
are linked. As explained in \cite{OSKh}, Lemma 3.6,
the linking number between two components is given by the number of twists in the unknot
between the attaching points of the two chords. In our case, because of the special position
of the unknot $U$ and the chords, this linking number is always zero or $\pm 1$, and can be easily
determined by untwisting the unknot (we can also pinpoint the sign of the linking number if we orient the surgery link by using 
the blackboard framing of $U$, as described in \cite{OSKh}). It is convenient to picture $U$ as a round circle and to
mark on it the attaching points for the chords given by the crossings. In general,  this algorithm does not uniquely determine the contact surgery diagram for  $(\Sigma(L),\xi_L)$, as the positioning of Legendrian unknots with 
respect to one another remains unclear even if all linking numbers are known. Indeed, the smooth link 
formed by these unknots does not determine the corresponding Legendrian link.  A precise algorithm can be obtained 
with a little more effort and is contained in \cite{HKP}. However, the present construction  provides a lot of 
useful information about the contact structure $\xi_L$.

\begin{example} \label{ex1}
Two examples are shown on Fig.~\ref{diags}: the double branched cover for the right-handed trefoil
(the braid $\sigma_1^3$) gives the (unique) Stein fillable structure on $-L(3,1)$,
and the double cover branched over the transverse unknot with $sl=-3$ is an overtwisted sphere. 
(It is not hard to see that the method from the preceding paragraph determines the surgery contact 
diagram uniquely for these simple cases. Alternatively, one can use the algorithm from \cite{HKP}.)
For the surgery diagrams here and later on, we choose the more familiar
contact form $\xi_{std}=\ker(dz-xdy)$, which is isotopic, but not identical to
the rotationally symmetric contact form used for the braid representation of
the transverse links. We hope that this will lead to no confusion.
\begin{figure}[ht]
\includegraphics[scale=0.70]{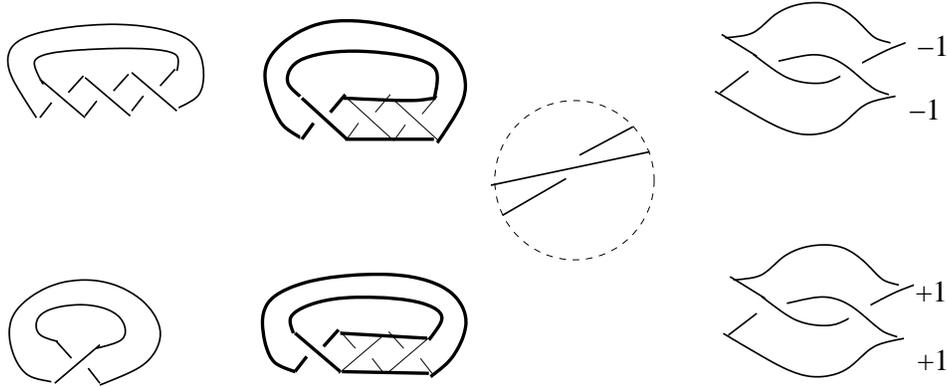}
\caption{Constructing surgery diagrams: the trefoil and the unknot.}
\label{diags}
\end{figure}
\end{example}

\begin{remark}
It is not necessary to single out the braid word that gives the unknot: we can as well
start from the trivial braid and obtain $(\Sigma(L),\xi_L)$ as a result of surgery
on $(\# S^1\times S^2, \xi_0)$ (where $\xi_0$ is the unique Stein fillable contact structure).
However, we find surgeries on the sphere more practical, especially for the next subsection.
\end{remark}

\subsection{Birman--Menasco braids}
We can now use the strategy from the previous subsection to construct surgery diagrams for the Birman--Menasco braids. We
make Theorem~\ref{bm1} more precise:

\begin{figure}[ht]
\vspace{.5cm}
\includegraphics[scale=0.72]{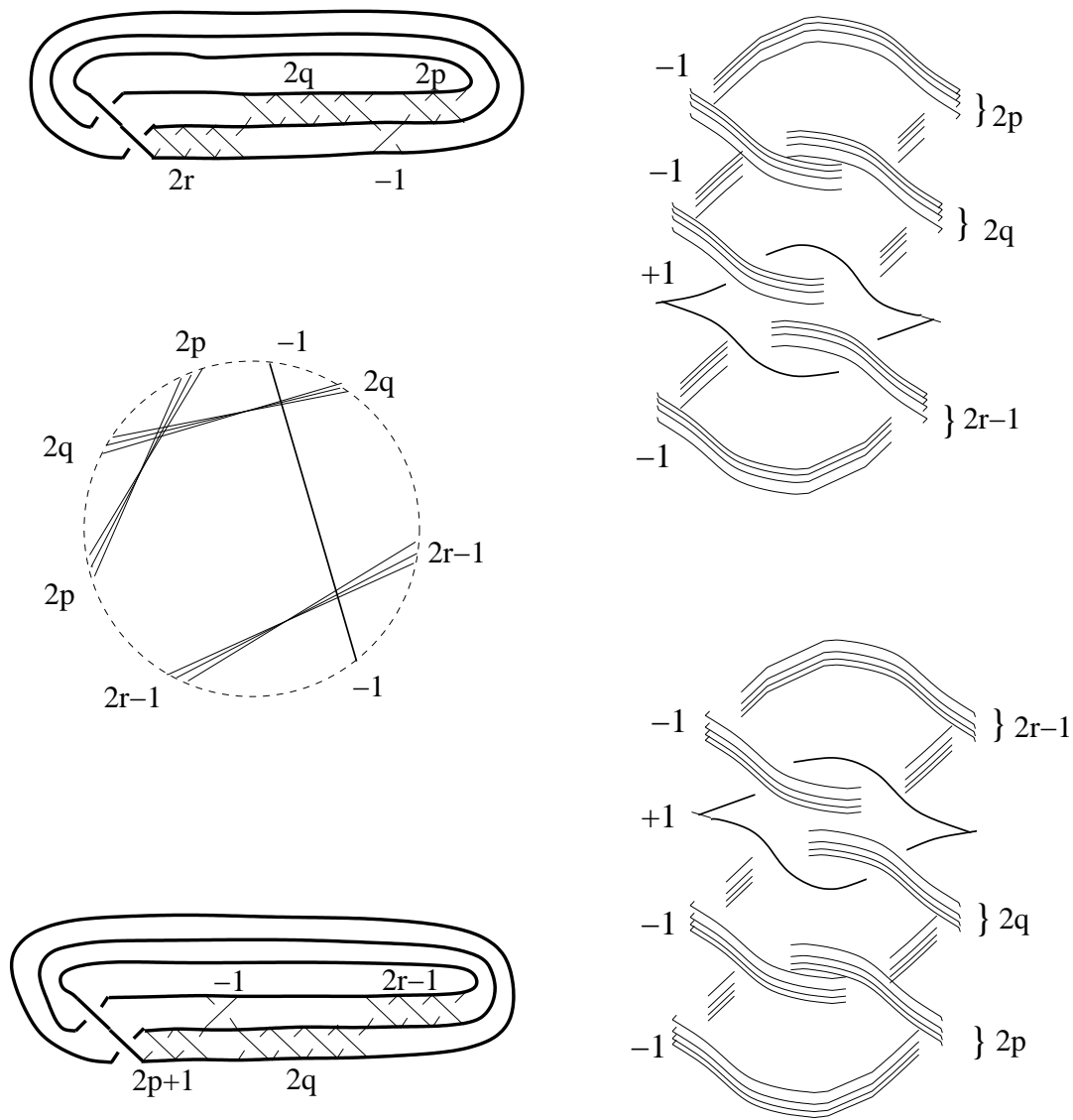}
\caption{Surgery diagrams for branched double covers of Birman--Menasco knots.} \label{birmen}
\end{figure}

\begin{theorem} \label{bmbraids} Let $p, q, r>1$ be integers with $p+1\neq q \neq r$,
and consider the transverse braids
$$
K_1=\sigma_1^{2p+1}\sigma_2^{2r}\sigma_1^{2q}\sigma_2^{-1}  \mbox{ and }  K_2=\sigma_1^{2p+1}\sigma_2^{-1}\sigma_1^{2q}\sigma_2^{2r}.
$$
The branched double covers  $(\Sigma(K_1),\xi_{K_1})$ and $(\Sigma(K_2),\xi_{K_2})$  are
contactomorphic.  
\end{theorem}

\begin{remark}
As mentioned in the introduction, $K_1$ and $K_2$ have the same self-linking number
and are isotopic as smooth knots, but not as transverse knots.
\end{remark}

\begin{proof} The two closed braids  representing $K_1$ and $K_2$ are isotopic (as closed braids) 
to 
$$
K_1=\sigma_1 \sigma_2^{2r} \sigma_1^{2q}\sigma_2^{-1} \sigma_1^{2p}  \mbox{ and } 
K_2=  \sigma_1 \sigma_2^{2p+1} \sigma_1^{-1}\sigma_2^{2q}\sigma_1^{2r-1}.
$$
Therefore, the two branched double covers are given by contact surgery on the same collection 
of Legendrian unknots, linked in the same way according to chord diagram on Figure~\ref{birmen}.  
This does not yet tell us what the surgery diagrams are; to determine the Legendrian links on which 
the surgery is performed, we use the open books techniques from \cite{HKP}. It turns out that 
$(\Sigma(K_1),\xi_{K_1})$ and $(\Sigma(K_2),\xi_{K_2})$ are  given by contact surgery 
on Legendrian links that are Legendrian mirrors of one another (see Figure~\ref{birmen}), and 
thus represent conjugate contact structures. These conjugate contact structures are in fact 
contactomorphic; we refer to \cite{HKP} for a proof.  

\end{proof}

\begin{remark}
We can also study the  manifold $(\Sigma(L),\xi_L)$ in terms of its open book
decomposition. Indeed, the page of such open
book is given by a branched double cover of a disk which is transverse to all
strings of the braid. The monodromy comes from the braid monodromy: each generator
$\sigma_i$ corresponds to a Dehn twist.  The manifold  $(\Sigma(L),\xi_L)$ can then
be exhibited as the boundary of an achiral Lefschetz fibration.
The open book description, together with the results of \cite{Gi}, would lead to an
alternative quick proofs for the next section. The surgery diagrams are more convenient
for our purposes because they fit nicely with the arguments of \cite{OSKh}.
\end{remark}

\section{Quasipositive braids and Stabilizations} \label{b-cover}
To prove Propositions \ref{qua-pos} and \ref{stab}, we turn attention to transverse links represented as quasipositive
braids, as well as to those which can be obtained as transverse stabilizations.

Recall \cite{Ru}  that a braid is called {\em quasipositive} if its braid
word is a product of conjugates of the form $w \sigma_i w^{-1}$,
where $w$ is an arbitrary element of the braid group.

\begin{proof}[Proof of Proposition \ref{qua-pos}]
Resolving a few positive crossings, we convert the braid representing $L$ into a braid
equivalent to a trivial one (of the same braid index). For the trivial braid $\OO$,
the contact manifold $(\Sigma(\OO), \xi_\OO)$  is a connected sum of several copies of
$(S^1\times S^2, \xi_0)$, which is Stein fillable. By Theorem \ref{resolve},  $(\Sigma(L),\xi_L)$
is obtained from $(\Sigma(\OO), \xi_\OO)$ by a sequence of Legendrian surgeries,
so it will be Stein fillable, too.

\end{proof}

Conversely, a stabilized transverse link gives some $(+1)$ surgeries in the contact surgery diagram for the double branched cover.


\begin{proof}[Proof of Proposition \ref{stab}] Suppose that a transverse link $L_{stab}$ is obtained as a transverse stabilization
of a link $L$. Then $L_{stab}$ can be represented
as a braid with a negative kink. Now use the algorithm  from section \ref{surgery}
to translate the braid representation into a contact surgery diagram. We see that
the diagram breaks into two pieces: two (+1) surgeries on two linked Legendrian unknots (cf. Example \ref{ex1}, Fig. \ref{diags}), isolated from everything else,
and
the surgery diagram for $(\Sigma(L), \xi_{L})$.
It follows that   $(\Sigma(L_{stab}), \xi_{stab})$ is the connected sum of  $(\Sigma(L), \xi_{L})$ and an overtwisted $S^3$, so it is overtwisted.
\end{proof}

\begin{cor} Let $T$ be  a transverse link smoothly isotopic to
a $(p, q)$ torus link, $p, q>0$. If $sl(T)=sl_{max}=pq-p-q$, then
$(\Sigma(T), \xi_{T})$ is Stein fillable. Otherwise $(\Sigma(T),
\xi_{T})$ is overtwisted.
\end{cor}

\begin{proof} We make use of the transverse simplicity of torus links \cite{Et1},
that is, the fact that
 the transverse $(p,q)$ torus link is uniquely determined by its
self-linking number. When $sl(T)=pq-p-q$, the link $T$ can be represented
as the obvious positive braid, so the branched double cover is
Stein fillable. For smaller values of $sl$, the transverse link $T$ is obtained as a
result of a few stabilizations of that braid, which proves
overtwistedness.
\end{proof}




\section{$\Spinc$ structures and the three-dimensional invariant}
\label{3dim}

A contact structure $\xi$ on a contact manifold $Y$ induces a
$\Spinc$ structure $\s$ on $Y$. When $c_1(\s)$ is torsion, there
is also the three-dimensional invariant $d_3(\xi)$ \cite{Go}. If
$(Y,\xi)$ is the boundary of an almost-complex 4-manifold $(X,J)$,
this invariant is given by
$$
d_3(\xi)=\f{1}{4}(c_1^2(J)-2\chi(X)-3\sign(X)).
$$

As described in section \ref{surgery}, given a braid presentation for
a transverse knot $K$, we can find a contact surgery description
of $(\Sigma(K), \xi_K)$. More precisely, assuming that the braid
has $b$ strings, and the braid word contains each $\sigma_i$ for
$i=1,\dots, b$  and  has $l+b$ positive entries and $m$ negative
entries, we see that $(\Sigma(K), \xi_K)$ is obtained from $(S^3,
\xi_{std})$ by a sequence of $l$ Legendrian surgeries and $m$ (+1)
contact surgeries, so that all the surgeries are performed on the
standard Legendrian unknot.

The results of \cite{DGS} can be used to understand the $\Spinc$
structure and the $d_3$ invariant for the contact manifold
$(\Sigma(L), \xi_L)$ for a transverse link $L$.

Let $X$ be the four-manifold obtained from $D^4$ by attaching the
2-handles as dictated by the $(\pm 1)$-surgery diagram. (If the
link $L$ has no negative crossings, $X$ is Stein.)
 Following
\cite{DGS}, consider an almost-complex structure $J$ defined on
$X$ in the complement of $m$ balls lying in the interior of the
(+1)-handles of $X$. As shown in \cite{DGS}, $J$ induces a
$\Spinc$ structure $s_J$ which extends to all of $X$. The $d_3$
invariant of $\xi_L$ can be computed as
\begin{equation}\label{d3}
d_3(\xi_K)=\f{1}{4}(c_1^2(\s_J)-2\chi(X)-3\sign(X))+m.
\end{equation}
This formula is very similar to the case where $(X, J)$ is
almost-complex, except that there is a correction term of $+1$ for
each (+1)-surgery.

Now, suppose that a $2$-handle is attached to the four-manifold
$X$  in the process of Legendrian surgery on a knot $K$, and
denote by $[S]$ the homology class that arises from the Seifert
surface of $K$ capped off inside the handle. It is well-known
\cite{Go} that $c_1(\s_J)$ evaluates on  $[S]$ as the rotation
number of the Legendrian knot $K$. Furthermore, it is shown in
\cite{DGS} that the same result is true for $(+1)$-contact
surgeries (for the $\Spinc$ structure $s_J$ on $X$ described
above). Thus, we see that in our case
 $c_1(\s_J)$ evaluates as
$0$ on each homology generator corresponding to either a $(-1)$ or
a $(+1)$ surgery, because the rotation number of the standard
Legendrian unknot is $0$. We conclude that $c_1(\s_J)=0$. Since
 $\s_\xi$ is the restriction of $\s_J$, we have
\begin{lemma} For any transverse link $L$, $c_1(\s_{\xi_L})=0.$
\end{lemma}

The Euler characteristic of the manifold $X$ is $1+\# \mbox{(2-handles)}$,
 which is $1+l+m$;  formula (\ref{d3}) simplifies as
$
d_3(\xi_L)=-\f{3}{4}\sigma(X)-\f{1}{2}(l-m)-\f{1}{2}=-\f{3}{4}\sigma(X)-\f{1}{2}sl(L)-\f{1}{2}.
$

Our next task is show that this expression depends on the
topological type of $L$ and $sl(L)$ only.

\begin{lemma} Suppose that two closed braids $L$ and $L'$ are
isotopic as smooth knots, and that $sl(L)=sl(L')$. Then
$d_3(\xi_L)= d_3(\xi_{L'})$.
\end{lemma}

\begin{proof} Since the braids $L$ and $L'$ give rise to isotopic
knots, by the classical Markov theorem \cite{Bi} $L'$  can be
obtained from $L$ by a sequence of braid isotopies and (positive
and negative) braid stabilizations and destabilizations. Braid
isotopies do not change the transverse link type, and neither do
positive stabilizations; this means that both the self-linking
number and the $d_3$ invariant remain unchanged. A negative
stabilization changes both $sl(L)$ and $d_3(\xi_L)$, but it is
easy to keep track of the changes. Indeed, the self-linking number
decreases by $2$. The  branched double cover $(\Sigma(L_{stab}),
\xi_{L_{stab}})$ is obtained from $(\Sigma(L), \xi_{L})$ by the
connected sum with an overtwisted sphere which is a (+1) contact
surgery on two linked Legendrian unknots (cf. Proposition \ref{stab}). It
follows from (\ref{d3}) that the $d_3$ invariant increases by $1$.
Since $sl(L)=sl(L')$, every negative stabilization must be
compensated by a destabilization, so that we have the same number
of stabilizations and destabilizations. Then, we must have
$d_3(\xi_L)= d_3(\xi_{L'})$.
\end{proof}

\section{Ozsv\'ath--Szab\'o invariants} \label{OSi}
In this section, we study the Ozsv\'ath--Szab\'o invariants of the contact
structures on the branched double covers. As the author learned
upon completion of this paper, the same question was independently studied by John Etnyre,
who obtained similar results.

\subsection{A brief review}
We quickly recall a few facts about the Heegaard Floer homology
groups here, referring the reader to the papers of Ozsv\'ath and
Szab\'o for details. We use coefficients in $\Z$.  Given a
$3$-manifold $Y$ equipped with a $\Spinc$ structure $\s$,
Ozsv\'ath and Szab\'o use a Floer-theoretic construction to define
a homology group $\h{HF}(Y, \s)$. Together, these groups form
$\widehat{HF}(Y)= \bigoplus_{\s\in\Spinc(Y)}\widehat{HF}(Y,\s)$.
Cobordisms induce a map on Floer homology. More precisely, a
$\Spinc$ cobordism $(W, \s)$ gives a map
$$
F_{W, \s}:\widehat{HF}(Y_1,\s|Y_1) \to \widehat{HF}(Y_2,\s|Y_2).
$$
The map $F_W:\h{HF}(Y_1)\to\h{HF}(Y_2)$ is defined by summing over
all $\Spinc$ structures.

 For a manifold $Y$ equipped with a
contact structure $\xi$, Ozsv\'ath and Szab\'o \cite{ContOS}
define a contact invariant $c(\xi)\in \widehat{HF}(-Y)$. The
element $c(\xi)$ is well-defined up to a sign, and lives  in
the $\Spinc$ component $\h{HF}(-Y, \s_{\xi}) $ associated to the
contact structure. Cobordisms given by Legendrian surgeries
respect $c(\xi)$.

 \begin{prop}\label{lcob}{\em(\cite{LS1})} Let $(Y', \xi')$ be obtained
from $(Y, \xi)$  by  a Legendrian surgery. Denote by  $W$  the
surgery cobordism, and  let $F_W:\h{HF}(-Y')\to \h{HF}(-Y)$ be the
associated map. Then $F_W(c(\xi'))=c(\xi)$.
\end{prop}

When $\s$ is a torsion $\Spinc$ structure, $\h{HF}(Y)$ is a
$\Z$-graded group. (Strictly speaking, this grading $\gr$ takes
values in $\Q$; it is a $\Z$-grading shifted by a rational
constant.) The degree of $c(\xi)$ is closely related to the
3-dimensional invariant of $\xi$.

\begin{prop} {\em(\cite{ContOS})} Suppose $\s_{\xi}$ is torsion. Then $c(\xi)$ is a
homogeneous element of degree  $\gr(c(\xi))=d_3(\xi)+\f{1}{2}$.
\end{prop}

Finally, we have the following important fact.
\begin{theorem}{\em(\cite{ContOS})}
{\em (1)} If the contact structure $\xi$ is overtwisted, then
$c(\xi)=0$.

{\em (2)} If $\xi$ is Stein fillable, then $c(\xi)\neq 0$. Indeed,
$c(\xi)\in \h{HF}(-Y)$ is a primitive element.
\end{theorem}

\subsection{Contact invariants of double covers} Combining the results of the previous sections
 with the properties of the Heegaard Floer contact invariant, we
immediately get the following propositions.

\begin{prop} Suppose that the transverse link $L$ is the result of transverse stabilization of another link.
Then $c(\xi_L)=0$.
\end{prop}

\begin{prop} Let the transverse link $L$ be represented as a quasipositive braid. Then $c(\xi_L)\neq 0$.
Indeed, it is a homogeneous primitive element in
$\h{HF}(-\Sigma(L))$, whose grading is given by
$$
\gr(\xi_L)=-\f{3}{4}\sigma(X)-\f{1}{2}sl(L),
$$
where $X$ is the 4-manifold described in section \ref{3dim}. The
latter expression depends on the topological type of $L$ and its
self-linking number only.
\end{prop}

\begin{prop} \label{cxito} Let the link $L$ be obtained from the link $L_+$ by resolving a positive crossing.
Let $W$ be the associated surgery cobordism between  the two branched double covers,
and  $F_W: \h{HF}(-\Sigma(L_+))\to \h{HF}(-\Sigma(L))$ the induced map on homology.
Then $F_W(c(\xi_{L_+}))=c(\xi_{L})$.
\end{prop}
We also have
\begin{prop} For the connected sum $L_1\# L_2$ of transverse links $L_1$ and $L_2$,
$c(\xi_{L_1\#L_2})=c(\xi_{L_1})\otimes c(\xi_{L_2})$.
\end{prop}
\begin{proof} The contact manifold $(\Sigma(L_1\#L_2), \xi_{L_1\#L_2})$ is the connected
sum of  $(\Sigma(L_1), \xi_{L_1})$ and $(\Sigma(L_2), \xi_{L_2})$,
so the proposition follows from the connected sum statement in
\cite{ContOS}.
\end{proof}

The following proposition is useful in calculations  and shows
that $c(\xi_L)=0$ in many cases.

\begin{prop} \label{neg=0} Suppose that the transverse link $L$ is represented by a
closed braid such that its braid word contains a factor of
$\sigma_i^{-1}$ but no $\sigma_i$'s for some $i>0$. (This means
that all the crossings in the braid diagram on the level between
$(i-1)$-th and $i$-th string are negative.) Then $c(\xi_L)=0$.
\end{prop}

\begin{proof} First of all, we delete all $\sigma_i^{-1}$ but one
from the braid word, obtaining a link that decomposes as a
connected sum of two links (connected by a negative crossing, the
$\sigma_i^{-1}$ that remains). Deleting further negative crossings
and inserting new positive crossings into both components of the
connected sum, we obtain a link $L'$ given by two positive torus
knots connected by a negative crossing. Topologically, the link
$L'$ is just the connected sum of two torus knots; as a transverse
link, it does not have maximal self-linking number (because we can
connect the two components by a positive crossing instead of a
negative one to increase $sl$). Connected sums of torus knots are
transversely simple \cite{EH2}, so $L'$ is the transverse
stabilization of another link. By Proposition \ref{stab},
$c(\xi_{L'})=0$. Repeated use of Proposition \ref{cxito} now
implies that $c(\xi_L)=0$.
\end{proof}

\begin{cor} Let $L$ be a transverse representative of a  negative torus
link. Then $c(\xi_L)=0$.
\end{cor}

\section{Relation with Khovanov homology} Now we explore the
connection to the Khovanov homology mentioned in the introduction.
We consider the Khovanov homology with $\Z/2\Z$ coefficients.

\begin{figure}[ht]
\includegraphics[scale=0.7]{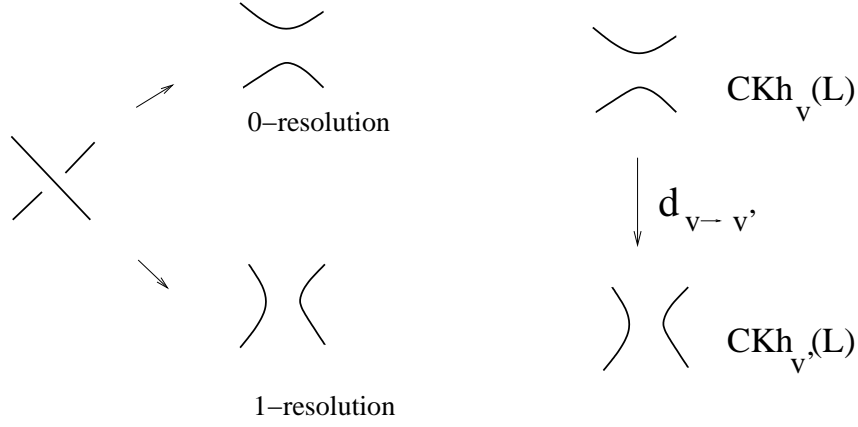}
\caption{Crossing resolutions and differentials.}
\label{crossings}
\end{figure}

We very briefly recall the relevant constructions. The starting
point for defining the Khovanov homology \cite{Kho} of a smooth
link $L\subset S^3$ is a link diagram (which we still denote $L$).
Let $n$ be the number of crossings in $L$.
Khovanov homology $Kh(L)$ is the homology of the chain complex 
$CKh(L)$, which is formed 
 by considering all ``complete resolutions" of $L$. Each
crossing can be resolved in two ways, (``0-resolution" and
``1-resolution", see Fig. \ref{crossings}), so that complete
resolutions can be written as $L_v$, indexed by $v\in [0,1]^n$.
 Each complete resolution
consists of a few disjoint circles; let $L_v$ consist of $m$
circles. To obtain the generators for $CKh(L_v)$, label each of
the circles by either $\um$ or $\up$ (that is, set
$CKh(L_v)=U^{\otimes m}$, where $U$ is a vector space with basis
$\{\um, \up\}$). The underlying space for the Khovanov complex is
then $CKh(L)=\oplus_{v\in \{0,1\}^n} CKh(L_v)$. The differential $d$
is defined as a sum of its components $d_{v\to v'}: CKh(v)\to
CKh(v')$ for all  $v$, $v'$ which are adjacent vertices of the
cube $[0,1]^n$, so that $L_v'$ can be obtained from $L_v$ by
changing a 0-resolution of one crossing into a 1-resolution.
The chain complex $CKh(L)$ is bi-graded; up to a correction term, 
the homological grading
of all elements in $CKh(L_v)$ is given by the number of zeroes among 
the coordinates of $v\in \{0,1\}^n$, so that $d$ raises the homological grading by 1. 
 We do
not describe the differential in detail, nor do we discuss
the quantum grading on $CKh(L)$. (The reader is referred to the original
paper \cite{Kho} or to surveys such as \cite{BN}.)

 The reduced complex $\widetilde{CKh}(L)$ is defined for a link with a
marked point. Each complete resolution now has one marked circle;
let $CKh_{\um}(L)$ be generated by those
$\uu_{\pm}\otimes\dots\otimes\uu_{\pm}$ that have the label $\um$
on the marked circle. Then, the reduced  complex is defined as
$\widetilde{CKh}(L)=CKh(L)/CKh_{\um}(L)$. The group $\Kh(L)$ is
the reduced Khovanov homology group.

Given a transverse link $L \subset S^3$, we define an invariant
$\psi(L)\in \Kh(L)$ \cite{Pla}. We represent $L$ as a transverse
braid and  take the oriented resolution $L_o$ of the braid diagram
(i.e. the resolution that consists of parallel strings). For unreduced Khovanov homology, 
we  pick the element
$
\um\otimes\um\otimes \dots\otimes \um \in CKh(L),
$
that is, we label every  component of the oriented resolution with a $\um$.  
For reduced homology, there is an isomorphism $ \widetilde{CKh}(L) \cong CKh_{\um}(L)$
obtained by replacing the label $\up$ on the marked component with a $\um$; we define
$\ps(L)$ to be the element of 
$\widetilde{CKh}(L)$ corresponding to  $ \um\otimes\ \dots\otimes \um \in CKh_{\um}(L)$ under this isomorphism. 
We show that $\ps(L)$ descends to
$\psi(L)\in\widetilde{Kh}(L)$ and gives an invariant of the
transverse link, well-defined up to a sign. The properties of
$\psi(L)$ are similar to those of $c(\xi_L)$; in fact, we prove
analogs of results of section \ref{OSi}. In the introduction, we
have conjectured that $\psi(L)$  ``corresponds" to $c(\xi_L)$ under
the spectral sequence of \cite{OSKh} which relates $\Kh(L)$ and
$\h{HF}(-\Sigma(L))$.

To be more precise, we consider the special case when the link $L$
is alternating; in this case, the spectral sequence collapses, 
 giving an isomorphism
 between $\Kh(L)$ and $\h{HF}(-\Sigma(L))$.
 We recall how this isomorphism
is established, describing $\h{HF}(-\Sigma(L))$ as the homology of
the following filtered chain complex \cite{OSKh}.
 Given a link diagram $L$, consider again all complete resolutions
 $L_v$, $v\in\{0,1\}^n$. For each $L_v$, let $Y_v=\Sigma(L_v)$ be
 the double cover of $S^3$ branched over $L_v$. (Since $L_v$ is simply
 the disjoint union of say $m$ circles, $Y_v=\#_{m} S^1\times S^2$.)
 Now, let $E^1=\oplus _{v\in\{0,1\}^n}\h{HF}(-Y_v)$ be the underlying space of the
 chain complex, and construct the differential $D=D^1$ as follows.
 As in Khovanov's theory, $D$ is the sum of its components
 $D_{v\to v'}:\h{HF}(-Y_v)\to \h{HF}(-Y_{v'})$ for all adjacent
 $v$, $v'$ (such that $L_v'$ is obtained from $L_v$ by
changing a 0-resolution of one crossing into a 1-resolution).
Then, $Y_v$ is obtained from  $Y_{v'}$ by a single 2-handle
attachment. The map $D_{v\to v'}$ is then defined as a map on
Heegaard Floer homology associated to the handle attachment
cobordism. The filtration grading on $E^1$ parallels the homological 
grading in Khovanov's theory; again,  on $\h{HF}(-Y_v)$ it is given by the number of zeroes
among the coordinates of $v\in\{0,1\}^n$ (it is convenient to introduce a correction term, too,
so that the two gradings are the same).
\begin{theorem} \cite{OSKh} \label{iso} 
Let $L$ be an alternating link, and fix its alternating diagram.
The homology of the filtered chain complex $(E^1, D)$ is
$\h{HF}(-\Sigma(L))$. On the other hand,  the associated graded complex of 
$(E^1, D)$ is isomorphic to $(\CKh(L), d)$. (For both theories, the coefficients are taken in $\Z/2\Z$.)  
\end{theorem}
The isomorphism between the two chain complexes
comes as part of the construction: indeed,
$\h{HF}(-Y_v)=\h{HF}(-S^1\times S^2)=\CKh(L_v)$, and the maps
$D_{v\to v'}$ are the same as $d_{v\to v'}$ under this
equivalence. When the link diagram is fixed, 
this provides a canonical isomorphism between $\Kh(L)$ and the associated 
graded group of $\h{HF}(-\Sigma(L))$  (The latter is non-canonically isomorphic 
to $\h{HF}(-\Sigma(L))$.)

\begin{remark} As notation suggests, $(E^1, D^1)$ is the first
term of a certain spectral sequence \cite{OSKh}. This spectral
sequence has ``higher order" differentials defined via maps
$D_{v\to v'}$, with $v$, $v'$ not necessarily adjacent. For a
general smooth link, its $E^2$ term gives $\Kh(L)$, and
$E^\infty=\h{HF}(-\Sigma(L))$.
\end{remark}

One of the key features of the Heegaard Floer theory is the
surgery exact triangle. In the correspondence of \cite{OSKh}, it
parallels the skein exact sequence for the Khovanov homology. The
skein exact sequence for a link $L$ relates $\Kh(L)$, $\Kh(L_0)$
and $\Kh(L_1)$, where $L_0$ and $L_1$ stand for the 0- and
1-resolution of a given crossing of $L$. The surgery exact
sequence for a three-manifold $Y$ and a framed knot $\gamma\in Y$
relates $\h{HF}(Y)$, $\h{HF}(Y_1(\gamma))$ and
$\h{HF}(Y_0(\gamma))$, where the manifolds $Y_1(\gamma)$ and
$Y_0(\gamma)$ are obtained as the result of 1- resp. 0-surgery on
$K$. Since crossing resolutions for links  induce surgeries on the
branched double covers, the surgery triangle relates the Heegaard
Floer homology groups of manifolds $\Sigma(L)$, $\Sigma(L_0)$ and
$\Sigma(L_1)$, and the resulting exact sequence looks very similar
to the skein sequence for the Khovanov homology. When the link $L$
is alternating  (and so are its resolutions $L_0$ and $L_1$), the
two exact sequences fit together nicely.

\begin{lemma} Let the link $L$ be given by an alternating diagram,
and let the links $L_0$ and  $L_1$ be obtained by  the 0- and
1-resolution of a given crossing. Then the diagram
\begin{equation}\label{diag-com}
\begin{CD}
 @>>> \Kh(L) @>>> \Kh(L_0) @>>> \Kh(L_1) @>>>  \\
&& @VVV @VVV @VVV \\
 @>>> \h{HF}(-\Sigma(L)) @>>> \h{HF}(-\Sigma(L_0)) @>>>
\h{HF}(-\Sigma(L_1)) @>>>
\end{CD}
\end{equation} commutes. (That is, the two squares shown commute,
and the third square implicit in the diagram also commutes.) Here
the maps on the Heegaard Floer homology are induced by handle
attachments, the maps on the Khovanov homology are induced by the
crossing change, and the vertical maps are the isomorphisms provided by
Theorem \ref{iso}; to simplify notation here and below, we write 
$\h{HF}(-\Sigma(L))$ for the corresponding associated graded group.
\end{lemma}

As Peter Ozsv\'ath explained to the author, the proof of this
Lemma follows from the arguments and techniques of  \cite{OSKh}.
(In fact, we should consider maps on the corresponding filtered 
and associated graded complexes.)
To avoid a lengthy review, we do not include this proof here.
Instead, we consider a simple example which illustrates the
commutativity of the diagram (\ref{diag-com}) and the interplay
between $\psi(L)$ and $c(\xi_L)$.

\begin{example} Consider the transverse unknots $L^+$ and $L^-$, given by braid diagrams
with one positive resp. one negative crossing, and their 0- and 1-resolutions.
The exact triangles in the Heegaard Floer and the Khovanov theory look as follows.
\begin{figure}[ht]
\includegraphics[scale=0.65]{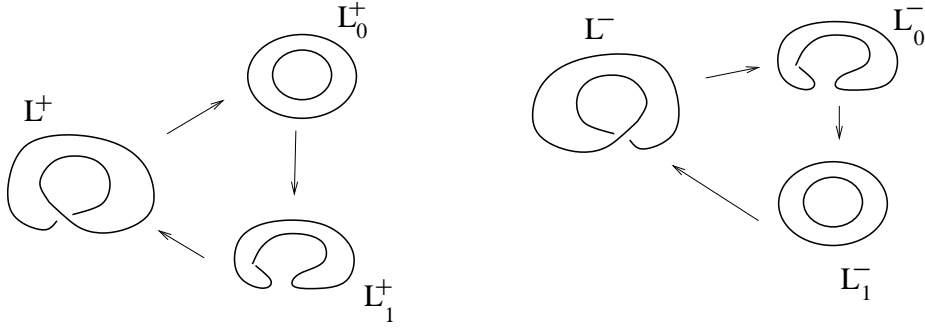}
\caption{The two unknots.}
\label{exa}
\end{figure}
For the link $L^+$, we have $\h{HF}(-\Sigma(L^+))=\h{HF}(S^3)=\Z/2\Z$, and also
$\h{HF}(-\Sigma(L^+_1))=\Z/2\Z$, while $\h{HF}(-\Sigma(L^+_0))=
\h{HF}(S^1\times S^2)=\Z/2\Z_{(-1/2)}\oplus \Z/2\Z_{(+1/2)}$, where the subscripts
indicate the gradings. The map $F:\h{HF}(-\Sigma(L^+))\to\h{HF}(-\Sigma(L^+_0))$
of the surgery exact sequence is the isomorphism between $\Z/2\Z$ and
$\Z/2\Z_{(+1/2)}\subset\h{HF}(-\Sigma(L^+_0))$;  the map $D:\h{HF}(-\Sigma(L^+_0))\to\h{HF}(-\Sigma(L^+_1))$ is an isomorphism
between  $\Z/2\Z_{(-1/2)}$ and $\Z/2\Z=\h{HF}(-\Sigma(L^+_1))$.
On the other hand, $\Kh(L^+)=\Z/2\Z$ is generated by a $\up$, and so is
$\Kh(L^+_1)=\Z/2\Z$. The group $\Kh(L^+_0)=\Z/2\Z \oplus \Z/2\Z$ is generated
by $\um\otimes\up$ and $\up\otimes\up$. The map $f:\Kh(L^+)\to\Kh(L^+_0)$
sends $\up$ to $\um\otimes\up$, while $d:\Kh(L^+_0)\to\Kh(L^+_1)$ sends $\up\otimes\up$
to $\up$. The groups  $\h{HF}(-\Sigma(L^+))$  and $\Kh(L^+)$ can be obtained as
homology of complexes $(\h{HF}(-\Sigma(L^+_0))\oplus \h{HF}(-\Sigma(L^+_1)), D)$
and $(\Kh(L^+_0)\oplus \Kh(L^+_1), d)$, respectively.

We observe:
\begin{enumerate}

\item $c(\xi_{L^+})=1\in\Z/2\Z=\h{HF}(-\Sigma(L^+))$, and $\psi(L^+)=\up\in \Z/2\Z=\Kh(L^+)$;

\item The maps $F$ and $f$ are associated to a positive crossing resolution,
so $F(c(\xi_{L^+}))=c(\xi_{L^+_0})=1\in\Z/2\Z_{+1/2}$,
$f(\psi(L^+))=\psi(L^+_0)=\um\otimes\up$;

\item $D(c(\xi_{L^+_0}))=0=d(\psi(L^+_0))$;

\item As an element of $\Kh(L^+_0)\oplus \Kh(L^+_1)=\CKh(L)$,
$\psi(L^+_0))=\um\otimes\up$ is the element $\ps(L^+)$ that represents the class
of $\psi(L^+)$. The identity $d(\ps(L^+))=0$ means that $\ps(L^+)$ is a cycle.
\end{enumerate}

For the link $L^-$, we have $\h{HF}(-\Sigma(L^-))=\h{HF}(S^3)=\Z/2\Z$, and now
$\h{HF}(-\Sigma(L^-_0))=\Z/2\Z$, while $\h{HF}(-\Sigma(L^-_1))=
\h{HF}(S^1\times S^2)=\Z/2\Z_{(-1/2)}\oplus \Z/2\Z_{(+1/2)}$. The map $F:\h{HF}(-\Sigma(L^-_1))\to\h{HF}(-\Sigma(L^-))$
of the surgery exact sequence is the isomorphism between $\Z/2\Z_{(-1/2)}\subset\h{HF}(-\Sigma(L^-_1))$ and $\Z/2\Z= \h{HF}(-\Sigma(L^-))$. Similarly, the map $D:\h{HF}(-\Sigma(L^-_0))\to\h{HF}(-\Sigma(L^-_1))$ is an isomorphism
between  $\Z/2\Z=\h{HF}(-\Sigma(L^-_0))$ and $\Z/2\Z_{(+1/2)}$.
The Khovanov homology  $\Kh(L^-)=\Z/2\Z$ is generated by a $\up$, and so is
$\Kh(L^-_0)=\Z/2\Z$, while $\Kh(L^-_1)=\Z/2\Z \oplus \Z/2\Z$ is generated
by $\um\otimes\up$ and $\up\otimes\up$. The map $f:\Kh(L^-_1)\to\Kh(L^-)$
sends $\up\otimes \up$ to $\up$, and $d:\Kh(L^-_0)\to\Kh(L^-_1)$ sends $\up$
to $\um\otimes \up$. In this case we have:
\begin{enumerate}

\item $c(\xi_{L^-})=0$, and $\psi(L^-)=0$;

\item The maps $F$ and $f$ are associated to a positive crossing resolution, so $F(c(\xi_{L^-_1}))=c(\xi_{L^-})=0$,
$f(\psi(L^-_1))=\psi(L^-)=0$;

\item The element $c(\xi_{L^-_1})=1\in\Z/2\Z_{+1/2}$ is the image under $D$
of the generator of $\h{HF}(-\Sigma(L^-_0)$; $\psi(L^-_1)$ is the image
of $\up\in\Kh(L^-_0)$ under $d$. As an element of $\CKh(L^-)$,
$\psi(L^-_1)$ is precisely $\ps(L^-)$, and the identity $\ps(L^-)=d(\up)$
means that the invariant $\psi(L^-)$ vanishes in $\Kh(L^-)$.

\end{enumerate}

\end{example}

\begin{proof}[Proof of Theorem \ref{psi=c}] Let the transverse link be represented
by a braid $L$ whose diagram is alternating. First of all, we
observe that such a braid enjoys little freedom. Indeed, suppose
there is a factor of $\sigma_1$ in the braid word. This gives a
positive crossing in the diagram, which means that the
``neighboring" crossing must be negative. Considering the crossings
one by one, we can conclude that all factors of $\sigma_1$ in the
braid word come with positive exponents, all factors of $\sigma_2$
come with negative exponents (that is, the braid word contains
$\sigma_2^{-1}$'s but no $\sigma_2$'s); all factors for $\sigma_3$
again have positive exponents, etc. Proposition \ref{neg=0} and
the analogous  proposition of \cite{Pla} imply that $c(\xi_L)=0$
and $\psi(L)=0$ unless the braid word for $L$ contains no
$\sigma_{2k}$ for any $k>0$. The latter situation means that $L$
is simply a disjoint union of $(2,n_i)$ torus links, and this is
the only case where the theorem needs proof.

When $L=T(2,n)$ is given by a 2-braid with $n$ positive crossings,
both $c(\xi_L)$ and $\psi(L)$ are non-zero. The map $f$ is associated 
to the resolution of a positive crossing; 
by definition,  $\psi(L)$ is a homogeneous element of homological degree 0 \cite{Pla}.
Recall that Khovanov homology  $\Kh(T(2,n))$  has a
special form \cite{Kho}, 
so that in every homological 
degree the component of $\Kh(T(2,n))$ is at most one-dimensional;
we can then say  that $\psi(L)$ is the (unique) element of the lowest homological grading 
in $\Kh(L)$.       Also, $\h{HF}(-\Sigma(L))=(\Z/2\Z)^n$ since $\Sigma(L)$ is the lens
space $-L(n,1)$. When  $n=1$, it follows that $c(\xi_L)=c_0(\xi_L)=\psi(L)$.


When  $n>1$, we proceed by induction. Consider the link
$L=T(2,n+1)$ and its two resolutions $L_0=T(2,n)$ and
$L_1=\unknot$.
\begin{figure}[ht]
\includegraphics[scale=0.7]{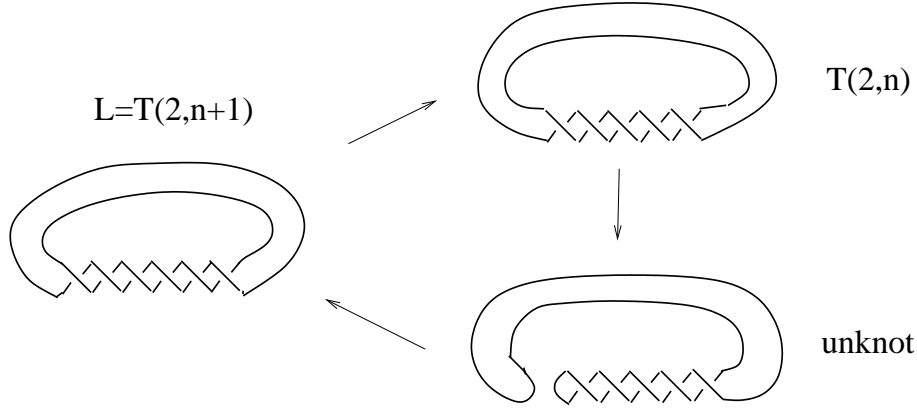}
\caption{The torus link $L$ and its resolutions.}
\label{toruslks}
\end{figure}

For these three links, (\ref{diag-com}) becomes
\begin{equation} \label{dia-tor}
\begin{CD}
\dots @>>> \Z/2\Z @>>> \Kh(L) @>f>> \Kh(T(2,n)) @>>> \dots \\
&& @VVV @VVV @VVV \\
\dots @>>> \Z/2\Z @>>> \h{HF}(-\Sigma(L)) @>F>>
\h{HF}(-\Sigma(T(2,n)) @>>>\dots
\end{CD}
\end{equation}
where $\Z/2\Z$ is the homology of the unknot. 

The image of $\Z/2\Z$ in $\Kh(L)$ lies in the homological degree higher than the minimal degree 0 
(where $\psi(L)$ lives). Then in the lower row of the above diagram, the image of $ \Z/2\Z$ 
in the associated graded group of $\h{HF}(-\Sigma(L))$  also lies in the higher filtration 
degree. The map $F$ on $\h{HF}(-\Sigma(L))$ is induced by the resolution of a positive crossing,
so $F(c(\xi_L))=c(\xi_{T(2,n)})$.  By the induction hypothesis, $c_0(\xi_{T(2,n)})$ lies in the lowest 
degree subquotient of $\h{HF}(-\Sigma(T(2,n)))$ and agrees with $\psi(T(2,n))$. The map $F$ preserves 
the filtration; it follows that it must send the lowest degree subquotient of $\h{HF}(-\Sigma(L))$ 
to the lowest  degree subquotient of $\h{HF}(-\Sigma(T(2,n)))$, since 
otherwise we couldn't have $F(c(\xi_L))=c(\xi_{T(2,n)})$. Then,
$c_0 (\xi_L)$ is also non-trivial, and so $c_0 (\xi_L)= \psi(L)$.



It remains to deal with the case when $L$ is a disjoint union of
torus links. Let $T_1= T(2, n_1)$, $T_2= T(2, n_2)$, and
$L=T_1\sqcup T_2$. Then, the contact manifold $(\Sigma(L), \xi_L)$
is the connected sum of $(\Sigma(T_1), \xi_{T_1})$, $(\Sigma(T_2),
\xi_{T_2})$ and $(S^1\times S^2, \xi_0)$. Then,
$\h{HF}(-\Sigma(L))=\h{HF}(-\Sigma(T(2, n_1)))\otimes
\h{HF}(-\Sigma (T(2, n_2)))\otimes 
\h{HF}(-S^1\times S^2)$, and $c(\xi_L)=c(\xi_{T_1})\otimes
c(\xi_{T_2}) \otimes c(\xi_0)$, where $c(\xi_0)\in
\h{HF}(-S^1\times S^2)=\Z/2\Z_{(-1/2)}\oplus \Z/2\Z_{(+1/2)}$ is
the generator of $\Z/2\Z_{(+1/2)}$.  Similarly, for the Khovanov
homology we have $\Kh(L)=\Kh(T_1)\otimes Kh(T_2)= \Kh(T_1)\otimes
\Kh(T_2)\otimes (\Z/2\Z\oplus \Z/2\Z)$, and
$\psi(L)=\psi(T_1)\otimes \psi(T_2)\otimes \um$. The
correspondence $\psi(L)=c_0(\xi_L)$ follows, since the isomorphism
between the homology groups maps the generator of
$\Z/2\Z_{(+1/2)}$ to the $\um$.
\end{proof}

\end{document}